\documentstyle[12pt]{article}

\newcommand{\auteur}{\noindent Nicole MESTRANO
{\footnotesize \newline Laboratoire Emile Picard
\newline UFR-MIG, Universit\'e Paul Sabatier\newline
31062 Toulouse CEDEX, France\newline
email: mestrano@picard.ups-tlse.fr}\vspace{.1cm}}

\newcommand{\numero}[1]{
\addtocounter{section}{1}
\begin{center}{\bf \thesection .\
#1\vspace{-.1in}}\end{center}
\setcounter{subsection}{0}
\setcounter{lemma}{0}\indent}

\newcommand{\subnumero}[1]{
\pagebreak[1]
\begin{center}{\em
#1}\nopagebreak\end{center} }

\newcommand{\eop}{\hfill $\Box$\vspace{.1in}}

\newtheorem{lemma}{Lemma}[section]

\newtheorem{proposition}[lemma]{Proposition}

\newtheorem{lemme}[lemma]{Lemme}
\newtheorem{theoreme}[lemma]{Th\'eor\`eme}
\newtheorem{corollaire}[lemma]{Corollaire}

\newtheorem{notations}[lemma]{Notations}
\newtheorem{remarque}[lemma]{Remarque}

\newtheorem{definition}[lemma]{D\'efinition}

\newcommand{\cc}{{\bf C}}

\newcommand{\zz}{{\bf Z}}
\newcommand{\pp}{{\bf P}}

\newcommand{\Ee}{{\cal E}}
\newcommand{\Ff}{{\cal F}}
\newcommand{\Gg}{{\cal G}}

\newcommand{\Oo}{{\cal O}}

\newcommand{\Mm}{{\cal M}}

\newcommand{\Jj}{{\cal J}}

\newcommand{\Cchi}{\begin{array}{c} \rule[-.12in]{0in}{.2in}
{\bf \chi}
\end{array}\!}

\begin{document}

\section*{Sur les espaces de modules des
fibr\'es vectoriels de rang deux sur des hypersurfaces
de $\pp ^3$}

\auteur

Soit $X$  une hypersurface lisse de degr\'e $\delta \geq 4$ dans $\pp ^3$
telle que $Pic(X)\cong \zz$.
On d\'esigne par $M_X(c_2)$ l'espace de modules des fibr\'es vectoriels
sur $X$ de rang $2$, de classes de
Chern $c_1 = 0 $ et $c_2$, semi-stables par rapport au diviseur hyperplan
(ou, ce qui est
\'equivalent, par rapport \`a n'importe quel  diviseur ample) \cite{Gieseker}.
Il est bien connu (Donaldson \cite{Donaldson}, Taubes \cite{Taubes},
Gieseker-Li \cite{GiLiJDG}, O'Grady \cite{OGrady}) que lorsque $c_2$
est suffisamment grand, l'espace  $M_X(c_2)$
est r\'eduit et irr\'eductible de dimension attendue.
Nous contribuons  ici \`a la recherche,  lorsque $c_2$ est petit,
de diff\'erentes composantes  irr\'eductibles.

On montre (corollaire \ref{L})
que pour tout entier $c_2 \geq \delta ^3/4 - \delta ^2/2 $
l'espace des modules $M_X(c_2)$ contient une composante
irr\'eductible r\'eduite de dimension attendue.
En utilisant un r\'esultat de O'Grady, on en d\'eduit
(corollaire \ref{M}) que pour tout entier $c_2$ tel que
$
\frac{1}{4}(\delta ^3-2\delta ^2)
\leq c_2 <
\frac{1}{3}(\delta ^3 - 9 \delta ^2 + 26 \delta - 3)
$
l'espace
$M_X(c_2)$ poss\`ede au moins deux composantes  irr\'eductibles
de dimensions distinctes.
(Pour $\delta \geq 27$, de tels entiers existent!)
On peut obtenir un r\'esultat analogue pour les fibr\'es vectoriels
avec $c_1 = 1$
(voir remarque \ref{R1}).

D'autre part, on prouve (th\'eor\`eme \ref{A1}) que pour $c_2 >
\frac{1}{12}(13 \delta
^3 - 24  \delta ^2 + 8 \delta)$,  le fibr\'e g\'en\'eral de la bonne
composante de
$M_X(c_2)$ que nous  construisons a la cohomologie naturelle.
\newline
Plus g\'en\'eralement, on d\'emontre (corollaire \ref{V}) que pour toute
surface
projective  lisse $X$ telle
que le fibr\'e canonique $K_X$ est de la forme $K_X = \Oo _X(k)$
pour un entier $k\in \zz$ o\`u $\Oo _X(1)$ est un fibr\'e ample,
si
$c_2$ est suffisamment grand alors, le fibr\'e g\'en\'eral de
l'unique composante de $M_X(c_2)$  a la
cohomologie naturelle o\`u $M_X(c_2)$ est l'espace de modules des
fibr\'es stables (par rapport \`a $\Oo _X(1)$) de classes de Chern
$c_1=0$ dans $Pic (X)$ et $c_2$.

Les fibr\'es vectoriels que nous construisons ici, et
qui sont dans une composante irr\'eductible r\'eduite de dimension attendue,
sont des  extensions
de $\Jj _{P/X}(\sigma )$ par $\Oo _X(-\sigma )$  o\`u $P$ est un   sous
sch\'ema de  dimension z\'ero en position sp\'eciale dans $X$.
On prend pour $P$ des intersections compl\`etes de $X$ avec des courbes
de $\pp ^3$  dont les id\'eaux ont les r\'esolutions les plus simples.
(Ces courbes ont \'et\'e  \'etudi\'ees
par de nombreux auteurs e.g.
\cite{Ellia}, \cite{Floystad}, \cite{Rao}, \cite{Walter}.)

C'est un grand plaisir de remercier Carlos Simpson
pour de nombreuses discussions et de lui d\'edier cet article.

\numero{Pr\'eliminaires}

Sauf mention  du contraire,
$X$  d\'esigne une hypersurface lisse de degr\'e $\delta \geq 4$ dans $\pp ^3$
telle que $Pic(X)\cong \zz$.

\begin{proposition}
\label{A}
Soient $\sigma $ un entier positif et  $C\subset \pp ^3$ une courbe lisse et
irr\'eductible de degr\'e $d$ telle que $H^1(C,\Oo _C(2\sigma -4))\neq 0$.
Quitte \`a translater $C$  par une hom\'eographie g\'en\'erale
de $\pp ^3$, si on  pose
$P:= C\cap X$ alors il existe  un fibr\'e vectoriel $E$ sur $X$
qui s'ins\'ere  dans la suite
exacte  $(\ast )$ suivante:
\vspace*{.5cm}
\newline
$(\ast )$
\vspace*{-1cm}
\newline
$$
0\rightarrow \Oo _{X}(-\sigma )
\rightarrow E \rightarrow \Jj _{P/X} (\sigma ) \rightarrow 0.
$$
On a $c_1(E)=0$ et $c_2(E)= \delta (d-\sigma ^2)$.
\end{proposition}
{\em D\'em:}
On sait (cf \cite{RHSRS}, Theorem 4.1) que la donn\'ee d'une courbe
$C\subset \pp ^3$
de degr\'e $d$, munie
d'une section non-nulle  $\eta \in H^0(C,\omega _C(4-2 \sigma ))$
est \'equivaute \`a celle
d'un faisceau reflexif $\Ff$  sur $\pp ^3$
de classes de Chern $c_1(\Ff)=2 \sigma$ et $c_2(\Ff)= d$
s'ins\'erant dans la suite exacte
$$
0 \rightarrow \Oo _{\pp ^3}
 \rightarrow \Ff \rightarrow \Jj _{C/\pp ^3}(2 \sigma ) \rightarrow 0.
$$
L'ensemble des points o\`u $\Ff$
n'est pas localement libre  est le sous sch\`ema
$Z$ de $C$ des z\'eros de la section $\eta$.
Soit $f : \pp ^3 \rightarrow \pp ^3$ une hom\'eographie
de $\pp ^3$ telle que $X$ ne rencontre pas $f^{-1}(Z)$.

On remplace $C$ par $f^{-1}(C)$ et on pose
$E:= f^{*}F(- \sigma)|_X$.
\eop

On \'etablit quelques propri\'et\'es de $E$.

\begin{lemme}
\label{B}
Soit $C\subset \pp ^3$ une courbe irr\'eductible lisse telle que
l'intersection $P := X\cap C$ soit transverse.
\newline
Si $H^0(\pp ^3 ,\Jj _{C/\pp ^3}(\tau ))= H^1(\pp ^3 ,\Jj _{C/\pp ^3}
(\tau -\delta ))=0$
alors $H^0(\pp ^3 ,\Jj _{P/X}(\tau ))=0$.
\end{lemme}
{\em D\'em.}
Soient $u\in H^0(\Jj _{P/X}(\tau ))$ et
$v \in H^0(\pp ^3, \Oo _{\pp ^3}(\tau ))$ un relev\'e de $u$
consid\'er\'e
comme une section de $\Oo _{X}(\tau )$.
La restriction $v|_C$ est une section de $\Oo _C(\tau )$ qui s'annulle sur
$P$
donc, la suite exacte
$$
0\rightarrow \Oo _{\pp ^3}(\tau -\delta )
\rightarrow \Oo _{\pp ^3}(\tau ) \rightarrow
\Oo _{X}(\tau ) \rightarrow 0
$$
restreinte \`a $C$ prouve qu'il
existe  un \'el\'ement  $w\in H^0(C,\Oo _C(\tau -\delta ))$ tel que $fw=v|_C$
o\`u   $f\in H^0(\pp ^3, \Oo _{\pp ^3}(\delta ))$ est l'\'equation de $X$.
Par hypoth\`ese
$H^1(\Jj _{C/\pp ^3}(\tau -\delta ))=0$, donc
$w$ s'\'etend en une section $w'\in H^0(\pp ^3, \Oo _{\pp ^3}(\tau -\delta ))$.
On a $(v-fw' )|_C = 0$. L'hypoth\`ese $H^0(\Jj _{C/\pp ^3}(\tau ))=0$ implique
$v-fw'=0$, et par cons\'equent $u = v|_X$=0.
\eop

\begin{corollaire}
\label{C}
Dans la construction la proposition \ref{A}, si on a
$\tau < \sigma$ et
\newline
$H^0(\pp ^3 ,\Jj _{C/\pp ^3}(\sigma + \tau ))=
H^1(\pp ^3 ,\Jj _{C/\pp ^3}(\sigma + \tau -\delta ))=0$,
alors  $H^0(X, E(\tau ))=0$.
\end{corollaire}
{\em D\'em.}
On applique le lemme pr\'ec\'edent en utilisant la suite exacte $(\ast )$.
\eop

\begin{corollaire}
\label{D}
Dans la construction de la proposition \ref{A},  si
\newline
$H^0(\pp ^3 ,\Jj _{C/\pp ^3}(\sigma ))=
H^1(\pp ^3 ,\Jj _{C/\pp ^3}(\sigma  -\delta ))=0$
alors  $E$ est stable.
\end{corollaire}
{\em D\'em.}
Puisque $E$ est un fibr\'e de rang $2$ sur une surface
dont le groupe de Picard est de rang un, il suffit de voir que
$H^0(X, E)=0$. pour cela, on applique
le corollaire pr\'ec\'edent avec $\tau = 0$.
\eop

\subnumero{Th\'eorie locale}

Nous rappelons ici quelques d\'efinitions et observations habituelles.

Soit  $[E]$ un point de l'espace des modules $ M_{X}^{st}(c_2)$ des
faisceaux sans
torsion stables de classes de Chern $c_1 =0$ et $c_2$,
repr\'esent\'e par le faisceau $E$.
Le morphisme de multiplication
$$
H^j(X, \Omega ^2_X) \rightarrow Ext^j(E,E\otimes \Omega ^2_X)
$$
donne, par dualit\'e de Serre (cf \cite{LePotier-Drezet} Proposition 1.2)
le morphisme {\em trace}
$$
Ext^i(E,E) \rightarrow H^i(X, \Oo _X).
$$
Comme le premier morphisme est injectif pour $j=0$, le morphisme trace
est surjectif pour
$i=2$. Le noyau $Ext^i(E,E)^o\subset Ext^i(E,E)$
du morphisme trace s'appelle la partie {\em sans
trace}.

L'espace tangent en $[E]$ \`a $M_X^{st}(c_2)$ est l'espace vectoriel
$Ext^1(E,E)^o$.
D'autre part,
le germe \'etale de l'espace des modules $M_X(c_2)$ en  $[E]$ est
isomorphe au germe \`a l'origine du sous-sch\'ema de $Ext^1(E,E)^o$
d\'efini par un
morphisme (nonlin\'eaire) $\Phi :Ext^1(E,E)^o \rightarrow Ext^2(E,E)^o$.
Ceci est
classiquement  bien connu. Une d\'emonstration purement alg\'ebrique est
indiqu\'ee par J. Le Potier dans ses notes \cite{LePotierLuminy}.

La {\em dimension attendue} de $M_X(c_2)$ en  $[E]$ est par
d\'efinition la diff\'erence
$$
{\rm dim\; att}  := {\rm dim} \; Ext^1(E,E)^o  - {\rm dim} \; Ext^2(E,E)^o
$$
Le th\'eor\`eme de Riemann-Roch  (avec l'hypothese $c_1=0$) fournit la valeur
$$
{\rm dim\; att} _E = 4c_2 - 3 \Cchi (\Oo _X).
$$

On a  $Ext^2(E,E)^o=0$ si et seulement si l'espace
$M_{X}(c_2)$ est lisse de dimension ${\rm dim\; att}$ au point $[E]$.
Dans ce cas, on dit que l'espace des modules est {\em bon} en $[E]$.

Une composante irr\'eductible non vide $Y\subset M_{X}(c_2)$
est {\em bonne} si l'espace des modules est bon au point g\'en\'erique de
$Y$ ce qui  \'equivaut \`a dire que $Y$ est r\'eduite de dimension
\'egale \`a
la dimension attendue.

Nous utiliserons la proposition bien connue suivante (cf  \cite{OGradyInvent}
ou \cite{GiLiJDG}).

\begin{proposition}
\label{E}
S'il existe une bonne composante  de $M_{X}(c_2)$, alors il
existe  une bonne composante de $M_{X}(c_2+1)$.
\end{proposition}
{\em D\'em.}
Soit $E'$ un  bon fibr\'e stable avec $c_1(E')=0$ et $c_2(E')=c_2$.
Soit $\Gg$ le faisceau gratte-ciel de rang $1$ au-dessus d'un
point $x\in X$. Soit $E_0$ le
noyau d'une surjection g\'en\'erale $E' \rightarrow \Gg$.
Alors $E_0$ est un  faisceau coh\'erent sans torsion
de classes de Chern  $c_1(E_0)=0$ et $c_2(E_0)=c_2 +1$.
Il a les m\^emes sous faisceaux que $E'$ donc il est stable.
\newline
Montrons que $E_0$ est bon c'est-\`a-dire que le morphisme trace
$t: Ext^2(E_0,E_0)  \rightarrow H^2(X,\Oo _X)$ est injectif.
\newline
En appliquant le foncteur $Ext(E', \ast )$ \`a
la suite exacte
$
0 \rightarrow E_0 \rightarrow E' \rightarrow \Gg \rightarrow 0
$
et en remarquant que $Ext^1( E', \Gg) = 0 $ (car $E'$ est
localement libre et $\Gg$ est un faisceau gratte ciel)
on obtient l'injection
$f: Ext^2( E',E_0) \rightarrow  Ext^2( E',E')$.
Par hypoth\`ese $E'$ est bon, donc le morphisme trace
$g: Ext^2(E',E')  \rightarrow H^2(X,\Oo _X)$ est injectif.
Soit maintenant
$h: Ext^2(E',E_0)  \rightarrow Ext^2(E_0,E_0)$
le morphisme  obtenu en  appliquant le foncteur $Ext( \ast ,E_0)$.
Comme $Ext^3( \Gg ,E_0) = 0 $, ce morphisme est surjectif.
On en conclut que $t$ est injectif en remarquant que
$t \circ h = g \circ f$.
\newline
Pour terminer la d\'emonstration de cette proposition,
il suffit maintenant de prouver que $E_0$ se d\'eforme en un
fibr\'e vectoriel $E$. En effet, par platitude, $E$ aura pour
classes de Chern  $c_1(E)=0$ et $c_2(E)=c_2 +1$
et il sera stable et
bon puisque ces propri\'et\'es sont des propri\'et\'es ouvertes.
\newline
On raisonne par l'absurde.
\newline
Supposons que pour toute d\'eformation
$\Ee$ de $E_0$ param\'etr\'ee par une courbe $T$
de  point g\'en\'eral $\eta$ et de point sp\'ecial
$s$ avec $\Ee_s = E_0$, le faisceau $\Ee_{\eta}$ ne
soit pas locallement libre.
Soit $\Ee$ une telle d\'eformation.
D\'esignons par $\Ee^{\ast \ast}$
le double dual du faisceau $\Ee$ et par $\Mm$ le
conoyau de l'injection de $\Ee$ dans $\Ee^{\ast \ast}$.
On a la suite exacte
$
0 \rightarrow \Ee \rightarrow \Ee^{\ast \ast} \rightarrow \Mm \rightarrow 0.
$
On a l'inclusion $E_0 \subset (\Ee^{\ast \ast})_s $ donc,
en remarquant que $E' = E_0 ^{\ast \ast}$, on voit que
$E' \subset  (\Ee^{\ast \ast})_{s} ^{\ast \ast}$. Or ce sont deux
faisceaux localement libre \'egaux hors codimension $2$,
d'o\`u $E' = (\Ee^{\ast \ast})_{s} ^{\ast \ast}$ et on a la double inclusion
$E_0 \subset (\Ee^{\ast \ast})_s \subset E'$.
Par hypoth\`ese  $longueur (E'/E_0) =1$,
donc $ longueur (E'/ (\Ee^{\ast \ast})_s) + longueur (\Mm_s) = 1$.
Par semicontinuit\'e,  $longueur (\Mm _{\eta}) \leq 1$ mais
on a suppos\'e que $\Ee_{\eta}$ n'est  pas locallement libre, donc
$longueur (\Mm _{\eta}) = 1$. Par semi-continuit\'e, on en d\'eduit
$longueur (\Mm_s) = 1$.
Par suite  $longueur (E'/ (\Ee^{\ast \ast})_s) = 0$
i.e. $E' = (\Ee^{\ast \ast})_s$ et $\Ee^{\ast \ast}$ est une d\'eformation de
$E'$.

Soit $U\subset M_X^{st}(c_2)$ un voisinage ouvert de $E'$ de dimension
$4c_2 - 3 \Cchi (\Oo _X)$.
Soit $Z _1\subset M_X^{st}(c_2 +1)$ l'ensemble des points correspondant
aux faisceaux
$F$ avec $longeur(F^{\ast\ast}/F)\geq 2$ et $Z_2\subset M_X^{st}(c_2+1)$
l'ensemble des
points correspondant aux faisceaux $F$ avec
\newline
$longeur (F^{\ast \ast}/F)=1$ et
$F^{\ast\ast} \not\in U$. Les sous-ensembles $Z_1$
et $Z_2$ sont constructibles.

L'argument ci-dessus montre que, pour $i \in \{ 1,2 \}$, la fermeture de
$Z_i$  ne contient pas $E_0$. En effet si $E_0$
\'etait dans la fermeture de $Z_i$  alors, il
existerait une d\'eformation
$\Ee$ de $E_0$ param\'etr\'ee par une courbe $T$
de point sp\'ecial
$s$ avec $\Ee_s = E_0$ et de  point g\'en\'eral $\eta$
avec   $\Ee_{\eta} \in Z_i$. D'apr\'es l'argument ci-dessus:
\newline
(i) \, $ longueur(\Ee^{\ast \ast}_{\eta}/\Ee_{\eta}) = 1$ donc
 $\Ee_{\eta} \not\in Z_1$ d'o\`u la contradiction pour $i=1$.
\newline
(ii) \, $\Ee^{\ast \ast}$ est une d\'eformation de
$E'$, donc  $\Ee^{\ast \ast}_{\eta} \in U$
et on ne peut pas avoir $\Ee_{\eta} \in Z_2$;
d'o\`u la contradiction pour $i=2$.

L'ensemble $Z_1\cup Z_2$  ne contenant pas $E_0$ dans sa fermeture,
son compl\'ementaire dans $M_X^{st}(c_2 +1)$,
contient un voisinage ouvert $V$ de $E_0$.
Quitte \`a restreindre $V$, puisque $E_0$ ne se d\'eforme pas en faisceau
localement
libre (par hypoth\`ese absurde), on peut supposer que les points de $V$
repr\'esentent des faisceaux $F$ non localement libres. Par suite,
si $F$ est repr\'esent\'e par un point de $V$, alors
$longeur(F^{\ast\ast}/F)=1$ et $F^{\ast\ast} \subset U$.
Les faisceaux $F$ s'ins\`erent donc dans des suites exactes de la forme
$$
0\rightarrow F \rightarrow F^{\ast\ast} \rightarrow M\rightarrow 0
$$
o\`u $M$ est un faisceau gratte-ciel de longeur $1$.
La dimension de l'ensemble de ces faisceaux $F$ est donc
inf\'erieure ou \'egale \`a
$ 3+ dim(U) = 3 + 4c_2 - 3\Cchi (\Oo _X)$.
Ceci nous fournit la contradiction cherch\'ee puisqu'on
sait que la
dimension de $M_X^{st}(c_2+1)$ en tout  point est au moins
\'egale \`a $4(c_2+1) - 3\Cchi (\Oo _X)$.
\eop

\begin{lemme}
Soit $E$ un fibr\'e vectoriel de rang $2$ et d\'eterminant trivial. On a
$$
Ext^i(E,E)^o = H^i(X, Sym ^2(E)).
$$
\end{lemme}
{\em D\'em.}
Comme $E$ est un fibr\'e on a $Ext ^i(E,E)^o= H^i (X,End ^o(E))$ o\`u $End
^o(E)$ est le noyau du morphisme trace $End(E)\rightarrow \Oo _X$.
D'autre part, $E$ \'etant de rang $2$ et de d\'eterminant trivial,
$E\cong E^{\ast}$, et par suite,
$End (E)= E^{\ast}\otimes E \cong E \otimes E$.  Via cet
isomorphisme le morphisme trace devient (au signe pr\`es)
$$
E\otimes E \rightarrow \bigwedge ^2E \cong \Oo _X.
$$
Le noyau du morphisme trace est donc le noyau du morphisme $E\otimes E
\rightarrow \bigwedge ^2E$, c'est-\`a-dire   $Sym ^2(E)$.
\eop

\begin{lemme}
Soit $E$ un fibr\'e s'ins\'erant dans la suite exacte  $(\ast )$.
Alors on a la suite exacte
\vspace*{.5cm}
\newline
$(\ast \ast )$
\vspace*{-1cm}
\newline
$$
0 \rightarrow E(-\sigma ) \rightarrow Sym ^2(E) \rightarrow
\Jj ^2_{P/X}(2\sigma )\rightarrow 0
$$
\end{lemme}
{\em D\'em.}
En tensorisant la suite exacte $(\ast )$ par $E$ on obtient le morphisme
$$
0 \rightarrow E(-\sigma ) \rightarrow E\otimes E
$$
qui, compos\'e avec la projection $E\otimes E\rightarrow Sym ^2(E)$
donne le morphisme
\newline
$f : E(-\sigma ) \rightarrow Sym ^2(E).$
C'est un exercice d'alg\`ebre lin\'eaire de voir que $f$
est injective au dessus du point g\'en\'erique, il s'ensuit que $f$
est une injection de faisceaux.

Soit $g: Sym ^2(E)
{\rightarrow} \Oo _X(2\sigma )$
le morphisme induit par la suite exacte $(\ast )$, en consid\`erant
$\Jj _{P/X}(\sigma ) $ comme un sous faisceau de $\Oo _X(\sigma )$.
On a  la suite suivante:
$$
0 \rightarrow E(-\sigma ) \stackrel{f}{\rightarrow} Sym ^2(E)
\stackrel{g}{\rightarrow} \Oo _X(2\sigma ).
$$
(C'est encore un exercice d'alg\`ebre lin\'eaire de voir que
la suite est exacte au dessus de  l'ouvert compl\'ementaire de $P$,
on en d\'eduit
l'exactitude sur $X$ en utilisant le fait que
toute section de $E(-\sigma )$ d\'efinie hors
codimension $2$ s'\'etend.)

L'image du morphisme naturel $h:
\Jj _{P/X}(\sigma )\otimes  \Jj _{P/X}(\sigma ) \rightarrow \Oo _X(2\sigma )
$
est le faisceau d'id\'eaux $\Jj ^2_{P/X}(2\sigma )$.
Donc, le diagramme commutatif suivant:
$$
\begin{array}{ccc}
E\otimes E & \rightarrow & \Jj _{P/X}(\sigma )\otimes  \Jj _{P/X}(\sigma ) \\
\downarrow & & \downarrow {}^{{}_h}\\
Sym ^2(E) & \stackrel{g}{\rightarrow} & \Oo _X(2\sigma )
\end{array}
$$
o\`u les fl\`eches verticale \`a gauche et   horizontale en haut sont
surjectives montre que l'image du morphisme $g$ est le faisceau
 $\Jj ^2_{P/X}(2\sigma )$.
\eop

\begin{proposition}
\label{F}
Soient $\sigma  $ un entier positif
et $P\subset X$ un sous-sch\'ema fini r\'eduit tels que:
\newline
i)\, \,\, $\delta -4 < 2\sigma$
\newline
ii)\, \, $H^0(X,\Jj _{P/X}(\delta - 4))= 0$, et
\newline
iii)\, $H^0(X, \Jj ^2_{P/X}(2\sigma + \delta - 4))=0.$
\newline
Alors tout fibr\'e vectoriel qui est extension de $\Jj _{P/X}(\sigma )$ par
$\Oo _X(-\sigma )$ est bon.
\end{proposition}
{\em D\'em.}
Soit $E$ un tel fibr\'e. D'apr\`es le lemme 1.6 et, par
dualit\'e il suffit de montrer que $H^0(Sym^2(E)\otimes
\Oo _X(\delta -4))=0$ (car
$E$ \'etant autodual, $Sym ^2(E)$ l'est aussi).
En tensorisant par $\Oo _X(\delta -4)$
la suite exacte $(\ast \ast )$ du lemme ci-dessus,
on obtient la suite exacte:
$$
0 \rightarrow E (\delta -4-\sigma ) \rightarrow Sym ^2(E) \otimes
\Oo _X(\delta -4)\rightarrow \Jj
_{P/X}^2(2\sigma + \delta - 4)\rightarrow 0.
$$
On est alors ramen\'es, gr\^ace \`a l'hypoth\`ese
$(iii)$, \`a prouver  que $H^0(E(\delta -4-\sigma ))=0$.

Pour ceci, on tensorise par $\Oo _X(\delta -4-\sigma )$
la suite exacte $( \ast )$.
Les hypoth\`ese $(i)$ et $(ii)$ permettent de conclure.
\eop

\subnumero{Construction via des points en position g\'en\'erale
(d'apr\`es O'Grady)}

\begin{proposition}
\label{G}
Soit $\delta \geq 14$ un entier.
Alors, pour tout entier $c_2$ avec
$$
\frac{1}{6}(\delta ^3 - 7 \delta )
 < c_2 <
\frac{1}{3}(\delta ^3 - 9 \delta ^2 + 26 \delta - 3),
$$
l'espace de modules $M_X(c_2)$
contient une composante irr\'eductible de dimension plus grande
que la dimension
attendue.
\end{proposition}
{\em D\'em.}
C'est un r\'esultat d'O'Grady (cf \cite{OGrady} Proposition 3.33)
appliqu\'e ici aux
hypersurfaces de $\pp ^3$. La condition $\delta \geq 14$
assure l' existence de tels entiers. La composante irr\'eductible obtenue ici
est celle qui contient les extensions  de $\Jj _{P/X}(\sigma )$ par
$\Oo _X(-\sigma )$
en prenant $ \sigma = 1$ et $P$ en position g\'en\'erale.
\eop

\numero{Construction de fibr\'es via certaines courbes de $\pp ^3$}

Soit $C$ une courbe localement intersection compl\`ete dans $\pp ^3$.
Rappelons les notations habituelles suivantes:
$$
s(C):= {\rm max} \{ s' \in \zz \; ; \;\; \forall s'' < s',\;\;
H^0 (\pp ^3, \Jj _{C/\pp
^3}(s'') =0\} ,  $$
et
$$
e(C) := {\rm max} \{ e'  \in \zz  \; ; \;\;  H^1 (C, \Oo _C(e') \neq 0 \} .
$$
Posons de plus:
$$
t(C):= {\rm max} \{ t'  \in \zz  \; ; \;\; \forall t'' < t',\;\;
H^1 (\pp ^3, \Jj _{C/\pp
^3}(t'') =0\}
$$

\begin{remarque}
\label{H}
Pour tout $e' \leq e(C)$ on a
$H^1 (C, \Oo _C(e')) \neq 0$. Notons aussi que
$t(C)$ est toujours au moins \'egal \`a $0$. (En fait, les
courbes  que l'on va utiliser auront $t(C)= \infty$.)
\end{remarque}

\begin{lemme}
\label{I}
Soit $C\subset \pp ^3$ une courbe irr\'eductible lisse telle que
l'intersection $P:= C\cap X$ soit transverse.  Si
$$
H^1(\pp ^3, \Jj _{C/\pp ^3}(n-\delta ))=
H^0(\pp ^3, \Jj _{C/\pp ^3}^2(n))=
H^0(C, N^{\ast}_{C/\pp ^3}(n-\delta ))=0
$$
alors
$$
H^0(X, \Jj _{P/X}^2(n))=0.
$$
\end{lemme}
{\em D\'emonstration}
Soit $f\in H^0(X, \Jj _{P/X}^2(n))$ consid\'er\'ee comme
section de $\Oo _X(n)$.
La suite exacte \vspace*{.5cm}
\newline
$(\ast \ast \ast )$
\vspace*{-1cm}
\newline
$$
0 \rightarrow \Oo _{\pp ^3}(n-\delta ) \rightarrow \Oo _{\pp ^3}(n)
\rightarrow \Oo _X(n) \rightarrow 0
$$
et le
fait que $H^1(\pp ^3, \Oo _{\pp ^3}(n-\delta ))=0$ impliquent
que $f$ s'\'etend en une section $g\in H^0(\pp ^3,
\Oo _{\pp ^3}(n))$.
\newline
La restriction $g|_C$ est une section de $\Oo _C(n )$ qui s'annule sur $P$
donc c'est une section de $\Oo _C(n-\delta )$
(cf. la suite exacte $
0 \rightarrow \Oo _{C}(n-\delta ) \rightarrow \Oo _{C}(n)
\rightarrow \Oo _P(n) \rightarrow 0
$
obtenue en restreignant $(\ast \ast \ast)$ \`a $C$).
\newline
Or la suite exacte
$$
0 \rightarrow \Jj _{C/ \pp ^3}(n-\delta ) \rightarrow \Oo _{\pp ^3}(n-\delta )
\rightarrow \Oo  _{C}(n-\delta ) \rightarrow 0
$$
et l'hypoth\`ese
$H^1(\pp ^3, \Jj _{C/\pp ^3}(n-\delta ))=0$
assurent l'existence d'une section
\newline
$h\in H^0(\pp ^3,
\Oo _{\pp ^3}(n-\delta )) $ telle que
$h|_C = g|_C$ (consid\'er\'ees comme
section de $\Oo _{C}(n-\delta )$).
\newline
Donc, quitte \`a remplacer $g$ par $g - h'$
o\`u $h'$ est l'image de $h$ par l'injection
\linebreak
$
0 \rightarrow \Oo _{\pp ^3}(n-\delta ) \rightarrow \Oo _{\pp ^3}(n)$,
on peut supposer que
$g|_C=0$ i.e.  $g \in H^0(X, \Jj _{C/\pp ^3}(n))$
(tout en gardant l'hypoth\`ese $g|_X=f$ i.e.  $g|_X$ s'annule
deux fois le long de $P = C \cap X$).
\newline
Soit $g'$  la d\'eriv\'e normale de $g$ le long
de $C$ c'est-\`a-dire l'image de $g$ dans \linebreak
$H^0(X, \Jj _{C/\pp ^3}(n) / \Jj _{C/\pp ^3}^2(n)) =
H^0(X, N^{\ast}_{C/\pp ^3}(n))$.
Comme $C$ et $X$ se coupent transversalement on a:
$$
\Jj _{P/X}/\Jj ^2_{P/X} = (\Jj _{C/\pp ^3}/\Jj ^2_{C/\pp ^3} )
\otimes _{\Oo _{\pp ^3}}
\Oo _X
= (\Jj _{C/\pp ^3}/\Jj ^2_{C/\pp ^3} )\otimes _{\Oo _C}
\Oo _P = N^{\ast}_{C/\pp ^3}|_P.
$$
\newline
Donc l'hypoth\`ese $g|_X=f \in H^0(X, \Jj _{P/X}^2(n))$ donne
$g'|_P = 0$.
\newline
Par suite,
$g' \in H^0(C, N^{\ast}_{C/\pp ^3}(n-\delta ))$.
\newline
Or on a suppos\'e $ H^0(C, N^{\ast}_{C/\pp ^3}(n-\delta ))=0$.
Donc $g'=0$, d'ou
$g\in  H^0(\pp ^3, \Jj _{C/\pp ^3}^2(n))$.
\newline
On a aussi suppos\'e $H^0(\pp ^3, \Jj _{C/\pp ^3}^2(n))=0$
donc $g=0$ d'o\`u $f=0$.
\eop

\begin{proposition}
\label{J}
Soit $\sigma >0$ et $C\subset \pp ^3$ une courbe irr\'eductible
lisse de degr\'e
$d$.  Quitte \`a translater $C$ par une hom\'eographie
g\'en\'erale de $\pp ^3$, si on pose $P:= C\cap X$et si
\newline
a)\, $2\sigma - 4 \leq e(C)$;
\newline
b)\, $\sigma < s(C)$ et $\sigma - \delta < t(C)$;
\newline
c)\, $\delta - 4 < 2 \sigma $;
\newline
d)\, $\delta - 4 < s(C) $;
\newline
e)\, $2\sigma - 4 \leq t(C)$;
\newline
f)\, $H^0(\pp ^3, \Jj^2_{C/\pp ^3}(2\sigma + \delta - 4))=0$; et
\newline
g)\, $H^0(C, N^{\ast}_{C/\pp ^3}(2\sigma - 4))=0$;
\newline
alors il existe un fibr\'e vectoriel $E$ extension de
$\Jj _{P/X}(\sigma )$ par
$\Oo _X(-\sigma )$. De plus $E$ est stable, bon et a pour classes de Chern
$c_1(E)=0$ et $c_2(E)= \delta (d-\sigma ^2)$.
\end{proposition}
{\em D\'emonstration}
L'existence du fibr\'e $E$ et le calcul de ses
classes de Chern se d\'eduisent de la proposition \ref{A},
gr\^ace \`a l'hypoth\`ese $(a)$. La stabilit\'e est cons\'equence
du corollaire \ref{D} via l'hypoth\`ese $(b)$.
Les hypoth\`eses $(e)$, $(f)$ et $(g)$ permettent d'utiliser le
lemme \ref{I} pour en d\'eduire
$H^0(X, \Jj _{P/X}^2(2\sigma + \delta - 4))=0$,
puis la proposition \ref{F} pour en d\'eduire (gr\^ace \`a $(c)$ et $(d)$)
que $E$ est {\em bon}.
\eop

\bigskip

Pour obtenir, pour de nombreux entiers $c_2$, une {\em bonne} composante
irr\'eductible de $M_X(c_2)$ de dimension attendue, il reste \`a
d\'eterminer  des courbes
$C$ qui nous donneront une grande famille d'entiers $\sigma$
v\'erifiant les conditions de la proposition 2.3.

\begin{proposition}
\label{K}
Pour tout entier $s \geq 1$, il existe
une courbe irr\'eductible lisse  $C$ de $\pp^3$ avec
 $s(C)= s$, $t(C)=\infty$ et $e(C) = s-3$
telle que
$H^0(C, N^{\ast}_{C/\pp ^3}(\tau ))=0$ pour $\tau < s$ et
$H^0(\pp^3, \Jj _{C/\pp ^3}^2(\tau ))=0$ pour $\tau < 2s$.
\end{proposition}

\noindent
{\em D\'emonstration:}
Les {\em courbes d\'eterminentielles} sont les courbes $C$ dont la
r\'esolution
de l'id\'eal est de la forme
$$
0\rightarrow \Oo _{\pp ^3}(-s-1) ^s \rightarrow \Oo _{\pp ^3}(-s)^{s+1}
\rightarrow \Jj _{C/\pp ^3} \rightarrow 0
$$
pour un entier $s$.
Floystad (cf \cite{Floystad}) et, ind\'ependamment,
Walter  (cf \cite{Walter}) ont
d\'emontr\'e  que
de telles courbes lisses existent pour tout $s\geq 1$. Il est  facile de voir
que  $s(C)= s$, $t(C)=\infty$ et $e(C) = s-3$ (cf
\cite{Floystad} ou \cite{Walter}). Rao
montre ( cf \cite{Rao} 1.12) d'une part que
$H^0(C, N^{\ast}_{C/\pp ^3}(\tau ))=0$ pour $\tau < s$, et, d'autre part,
que l'id\'eal $\Jj ^2_{C/\pp ^3}$ a la r\'esolution
$$
0\rightarrow \Oo _{\pp ^3}(-2s-2) ^a \rightarrow
\Oo _{\pp ^3}(-2s-1) ^b
\rightarrow \Oo _{\pp ^3}(-2s)^c \rightarrow \Jj ^2_{C/\pp ^3} \rightarrow 0.
$$
En utilisant le
fait que $H^i(\pp ^3, \Oo _{\pp ^3}(m))=0$ pour $i=1,2$,
on en d\'eduit facilement
que  $H^0(\pp^3, \Jj _{C/\pp ^3}^2(\tau ))=0$ pour $\tau < 2s$.
\eop

\begin{remarque}
\label{R20}{\bf :}
\end{remarque}
Dans une  premi\`ere version de cet article, nous avons donn\'e un argument
g\'eom\'etrique pour prouver les propri\'et\'es
$H^0(C, N^{\ast}_{C/\pp ^3}(\tau ))=0$ pour $\tau < s$ et
$H^0(\pp^3, \Jj _{C/\pp ^3}^2(\tau ))=0$ pour $\tau < 2s$.  Cet argument
\'etait bas\'e sur une courbe d\'eterminentielle sp\'eciale
singuli\`ere qui est une r\'eunion de droites (``stick figure'')
d\'ej\`a utilis\'ee par Fl{\o}ystad et Walter, avec un argument
combinatoire pour
les annulations.  Walter m'a ensuite fait remarquer que ces
propri\'et\'es ont
\'et\'e prouv\'ees par Rao.

\smallskip

\begin{corollaire}
\label{L}
Soit $X$ une hypersurface de degr\'e $\delta \geq 4$ dans $\pp^3$.
Pour tout entier $c$ avec $c\geq \delta ^3/4 - \delta ^2/2 $
il existe un bon fibr\'e stable sur $X$
avec $c_1=0$ et $c_2=c$.
(Si $\delta$ est impair le r\'esultat reste vrai pour
 $c\geq \delta ^3/4 - \delta ^2 + \frac{3}{4}\delta $.)
\end{corollaire}

\noindent
{\em D\'emonstration}
D'apr\`es la proposition \ref{E}, il suffit de construire un bon
fibr\'e stable
avec $c_1=0$ et $c_2 = \frac{1}{4}\delta ^2 (\delta - 2)$ si
$\delta$ est pair
(resp. $c_2 = \frac{1}{4}\delta (\delta -1)(\delta -3)$ si
$\delta$ est impair).
Pour cela posons
$
s= \delta -3 $  si $\delta $ est impair et $ s = \delta -2 $
si $\delta $ est pair.
Soit $\sigma = s/2$.
Soit $C$ une courbe comme dans la proposition \ref{K} On v\'erifie
facilement que les conditions de la proposition \ref{J} sont
satisfaites et on obtient
un bon fibr\'e stable $E$ ayant pour classes de Chern $c_1(E)=0$ et
$$
c_2(E)= (s(s+1)/2 - \sigma ^2)\delta
$$
$$
= (\sigma (2\sigma + 1) - \sigma ^2)\delta
$$
$$
= \delta \sigma (\sigma +1).
$$
En particulier, si $\delta$ est pair on a
$$
c_2(E)= \frac{1}{4}\delta ^2 (\delta - 2),
$$
et si $\delta $ est impair on a
$$
c_2(E) = \frac{1}{4}\delta (\delta -1)(\delta -3).
$$
\eop

\begin{corollaire}
\label{M}
Pour tout entier $c_2$
dans l'intervalle non vide suivant
$$
\frac{1}{4}(\delta ^3-2\delta ^2)
\leq c_2 <
\frac{1}{3}(\delta ^3 - 9 \delta ^2 + 26 \delta - 3)
\;\;\; \mbox {si} \; \; \delta \geq 28 \;\; \mbox{est un entier pair}
$$
ou
$$
\frac{1}{4}(\delta ^3-4\delta ^2 + 3\delta )
\leq c_2 <
\frac{1}{3}(\delta ^3 - 9 \delta ^2 + 26 \delta - 3)
\;\;\; \mbox {si} \;\; \delta \geq 21 \;\; \mbox{est un entier impair}
$$
l'espace des modules $M_X(c_2)$
contient une composante irr\'eductible
de dimension \'egale \`a la dimension
attendue et une autre
de dimension plus grande que la dimension
attendue.
\end{corollaire}
{\em D\'emonstration} C'est une cons\'equence directe
 de la proposition \ref{G} et du corollaire \ref{L}.
La condition $\delta \geq 28$ (resp.  $\delta \geq 21$)
assure l'existence d'entiers compris entre
$$
\frac{1}{4}(\delta ^3-2\delta ^2)
 \;\;\; \mbox{et} \;\;\;
\frac{1}{3}(\delta ^3 - 9 \delta ^2 + 26 \delta - 3)
$$
resp.
$$
\frac{1}{4}(\delta ^3-4\delta ^2 + 3\delta )
\;\;\; \mbox{et} \;\;\;
\frac{1}{3}(\delta ^3 - 9 \delta ^2 + 26 \delta - 3).
$$
\eop

\begin{remarque}
\label{R1}
{\bf :}
\end{remarque}
Nous avons utilis\'e la construction d'O'Grady dans le cas
o\`u elle donne des fibr\'es
{\em stables}.  On peut de m\^eme construire des extensions
de $\Jj _{P/X}$ par $\Oo _X$
(donc $H=0$  suivant les notations de \cite{OGrady}). On obtient
alors une  famille de
fibr\'es {\em semi-stables} de dimension plus grande que la
dimension attendue quand
$$ \frac{1}{6}(\delta ^3 - 6 \delta ^2 + 11 \delta ) < c_2 <
\frac{1}{3}(\delta ^3 - 6 \delta ^2 + 11 \delta -3).
$$
D'o\`u l'existence de deux composantes dans l'espace
des modules des faisceaux
sans torsion, l'une bonne et contenant
des fibr\'es stables, l'autre de
dimension plus grande que la dimension attendue contenant
des fibr\'es semistables (mais on
ne sait pas si elle contient des fibr\'es stables), pour
$$
\frac{\delta ^3}{4} - \frac{\delta ^2}{2} \leq c_2 <
\frac{1}{3}(\delta ^3 - 6 \delta ^2 + 11 \delta -3)
\;\;\;\;\; \mbox{si} \;\; \delta \;\; \mbox{est pair}
$$
et
$$
\frac{\delta ^3}{4} - \delta ^2 +\frac{3}{4}\delta  \leq c_2 <
\frac{1}{3}(\delta ^3 - 6 \delta ^2 + 11 \delta -3)
\;\;\;\;\; \mbox{si} \;\; \delta \;\; \mbox{est impair}.
$$
Ces intervalles sont non-vides pour $\delta \geq 16$
si $\delta $ est pair et pour $\delta \geq 9$
si $\delta$ est impair.

\begin{remarque}
\label{R2}{\bf :}
\end{remarque}
Pour simplifier la r\'edaction nous n'avons trait\'e que le cas
$c_1 = 0$ (c'est-\`a-dire, quitte \`a tensoriser par un fibr\'e en droites
ad\'equat, le cas $c_1$ pair). La d\'emonstration s'adapte au cas
$c_1 = 1$ (c'est-\`a-dire au cas $c_1$ impair) et on obtient
l'existence de deux composantes irr\'eductibles de
$M_X(1,c_2)$, l'une de dimension plus grande
que la dimension attendue et l'autre bonne (toutes  deux
contenant des fibr\'es
stables), pour
$$
\frac{1}{4}\delta (\delta - 1)^2 \leq c_2 <
\frac{1}{6}(2\delta ^3 - 15 \delta ^2 + 37 \delta -6)
\;\;\;\;\; \mbox{si} \;\; \delta \;\; \mbox{est impair}
$$
et
$$
\frac{1}{4}\delta (\delta - 2)^2  \leq c_2 <
\frac{1}{6}(2\delta ^3 - 15 \delta ^2 + 37 \delta -6)
\;\;\;\;\;  \mbox{si} \;\; \delta \;\; \mbox{est pair} .
$$
Ces intervalles sont non-vides pour $\delta \geq 21$ si $\delta $ est
impair et
pour $\delta \geq 14$
si $\delta$ est pair.

\numero{Fibr\'es \`a cohomologie naturelle}

Soit  $X$ une surface projective lisse telle que le fibr\'e
canonique $K_X$ soit de la forme $K_X = \Oo _X(k)$
pour un entier $k\in \zz$ o\`u $\Oo _X(1)$ est un fibr\'e ample.
Soit $M_X(c_2)$ l'espace des modules des fibr\'es stables
(par rapport \`a $\Oo _X(1)$)
avec $c_1=0$ dans $Pic (X)$ et $c_2$ donn\'e.

\begin{definition}
{\bf :}
\end{definition}
Un fibr\'e $E$ sur $X$ a la {\em cohomologie naturelle} (cf \cite{AHRH}) si
pour tout $n$ l'un au plus des groupes
$H^i(E(n))$ est non nul
pour $i=0,1,2$.
\newline
Plus g\'en\'eralement soit $\beta \geq k/2$ un entier. On dira qu'un
fibr\'e $E$ sur $X$ a la {\em cohomologie $\beta$-naturelle} si
pour tout $n \geq \beta$ l'un au plus des groupes
$H^i(E(n))$ est non nul pour $i=0,1,2$.

\begin{remarque}
{\bf :}
\end{remarque}
Pour tout fibr\'e $E$ avec $c_1(E) = 0$ dans $Pic(X)$, la dualit\'e de Serre
et la condition $K_X = \Oo _X(k)$, donnent l'\'egalit\'e
$h^2(X,E(n)) = h^0(X,E(k -n))$ et $h^1(X,E(n)) = h^1(X,E(k -n))$.
Il suffit donc d'\'etudier la cohomologie des $E(n)$
pour les entiers $n \geq \frac{k}{2}$.
Si de plus $E$ a  la {\em cohomologie naturelle} alors
pour $n \geq \frac{k}{2}$  on a $ h^2(E(n)) = 0$ et la cohomologie
des $E(n)$ se
d\'eduit de la  caract\'eristique d'Euler $\Cchi (E(n)) =
2\Cchi (\Oo _X (n)) -c_2$.
Enfin on peut
remarquer que si $k$ est pair et  $\Cchi (E(\frac{k}{2}))>0$ alors $E$
ne peut pas
avoir la cohomologie naturelle.

Nous \'etudions maintenant le comportement de la fonction de Hilbert de
$E$ quand on
applique la construction de proposition \ref{E} pour passer de $c_2$ \`a
$c_2+1$.

\smallskip

\begin{notations}
{\bf :}
\end{notations}
Soit  $c_2$ un entier  tel qu'il existe   un bon fibr\'e
vectoriel stable $E'$  sur $X$ avec $c_1(E')=0$ et  $c_2(E')=c_2$.
Soit $\Gg$ le faisceau gratte-ciel de rang $1$ au-dessus d'un point
$x\in X$. Soit $E_0$ le
noyau d'une surjection g\'en\'erale $E' \rightarrow \Gg$. Soit $E$
un fibr\'e stable et
bon qui est une d\'eformation g\'en\'erale de $E_0$.

\smallskip

\begin{proposition}
\label{T}
Soit $\beta \geq k/2$ un entier.
Si $E'$ a la {\em cohomologie $\beta$-naturelle}, alors $E$ aussi.
\end{proposition}
{\em D\'em.} Soit $n \geq \beta \geq k/2$ un entier.
Notons tout d'abord l'\'egalit\'e  $H^2(E'(n)) = 0$. En effet,
si $H^2(E'(n)) \neq 0$, alors, par dualit\'e de Serre,
$H^0(E'(k-n)) \neq 0$ donc,
$\forall m \geq k-n \;  , \; H^0(E'(m)) \neq 0$.
En particulier  $H^0(E'(n)) \neq 0$
et ceci est impossible si $E'$ a la cohomologie $\beta$-naturelle.
De la suite exacte
$$
0 \rightarrow E_0 \rightarrow E' \rightarrow \Gg \rightarrow 0
$$
et par semi continuit\'e, on en d\'eduit  $H^2(E(n)) = 0$.
\newline
De m\^eme, si  $H^0(E'(n)) = 0$ alors $H^0(E(n)) = 0$.
\newline
Enfin, si $H^0(E'(n)) \neq 0$, alors  d'une part $H^1(E'(n)) = 0$
et, d'autre part, pour une surjection g\'en\'erale le morphisme
$$
H^0(E'(n)) \rightarrow H^0(\Gg (n))= \cc
$$
est surjectif. Donc, $H^1(E(n)) = 0$.
\eop

\begin{notations}
{\bf :}
\end{notations}
On fixe maintenant un entier $c'_2$  tel qu'il existe
un bon fibr\'e vectoriel stable $E'$ sur $X$ avec $c_1(E')=0$ et
$c_2(E')= c'_2$. (C'est toujours possible, il suffit d'avoir
$c'_2$ assez grand.)
\newline
Soit $\beta $ un entier $\geq \frac{k}{2}$ tel que $E'$
a la cohomologie $\beta$-naturelle. (Un tel $\beta $
existe car $H^1(E'(n))=H^2(E'(n))=0$ pour $n$ assez grand
d'apr\'es les th\'eor\`emes d'annulations de Serre.)
\newline
Pour tout entier $c_2\geq c'_2$ d\'esignons par $M_X'(c_2)$ la composante
irr\'eductible de $M_X(c_2)$ obtenue \`a partir de $E'$ via la construction
d'augmentation de $c_2$ (it\'er\'ee $c_2-c'_2$ fois) de la proposition
\ref{E}.

\begin{proposition}
\label{U}
Soit $E \in M_X'(c_2)$ un fibr\'e ayant
la cohomologie $\beta$-naturelle.
Si $\Cchi (E(\beta ))<0$, alors $E$ a la cohomologie naturelle.
\end{proposition}
{\em D\'em.} Soit $n$ un entier. On doit
\'etudier la cohomologie de $E(n)$.
Par dualit\'e de Serre et puisque $E$ a
la  cohomologie $\beta$-naturelle, on peut supposer
$k/2 \leq n < \beta$.
\newline
L'hypoth\`ese
$\Cchi (E(\beta )) < 0 $ donne $H^1 (E(\beta )) \neq 0 $
et donc, $H^0 (E(\beta )) = 0 $. On en d\'eduit:
$$
\forall m \leq \beta \; , \; H^0 (E(m)) = 0.
$$
D'o\`u
$H^0 (E(n)) = 0$ et, par le m\^eme argument que pr\'ec\'edemment,
$H^2 (E(n)) = 0$.
\eop

\smallskip
\begin{corollaire}
\label{V}
On pose $\gamma := 2 \Cchi (\Oo _X(\beta ))$.
Si $c_2 > \gamma $ alors le fibr\'e g\'en\'eral $E$ de la composante
$M'_X(c_2)$ a la cohomologie naturelle.
En particulier quand $c_2\gg 0$ le fibr\'e g\'en\'eral de l'unique
composante de $M_X(c_2)$ a la cohomologie naturelle.
\end{corollaire}
{\em D\'em.}
D'apr\`es le choix de $\beta$ et  la proposition \ref{T},
le fibr\'e g\'en\'eral de
$M'_X(c_2)$ a la cohomologie $\beta$-naturelle.  L'hypoth\`ese $c_2(E)>
\gamma $
donne $\Cchi (E(\beta ))= 2\Cchi (\Oo _X(\beta ))-c_2 < 0$, donc $E$ a la
cohomologie naturelle  d'apr\`es la proposition \ref{U}.
\eop

\bigskip
Revenons au cas o\`u $X$ est une hypersurface lisse de degr\'e
$\delta$ dans $\pp^3$, avec $Pic (X)=\zz$. On a alors
$K_X = \Oo _X(\delta -4)$ si $\Oo _X(1)$ est le fibr\'e hyperplan.
La composante $M'_X(c_2)$ est celle que nous avons construite.
Les  propositions suivantes permettent de d\'eterminer $\gamma$.

\begin{proposition}
\label{R}
Soit $E'$ le fibr\'e vectoriel construit en $2.6$ avec,
suivant la parit\'e de $\delta$,
$c_2(E')=\delta ^3/4-\delta ^2/2$
ou $c_2(E')=\frac{\delta ^3}{4} -\delta ^2 +
\frac{3\delta }{4}$.
Soit $\beta $ le  plus petit entier $>  \frac{3}{2} \delta -4$.
Alors, $E'$ a la cohomologie $\beta$-naturelle.
\end{proposition}
{\em D\'em.}
Nous allons prouver que pour tout entier $n> \frac{3}{2}\delta -4$ on a
$h^1(E'(n))= h^2(E'(n))=0$.

Par construction et d'apr\`es le corollaire $1.3$, on a $h^0(E'(n))=0$ pour
$n< (\delta - 3)/2 $.  Par dualit\'e, on en d\'eduit
 $h^2(E'(n))=0$ pour $n> (\delta - 4)/2$.

D'autre part,  on sait que pour les courbes $C$
consid\'er\'ees ici, $H^1(\Jj _{C/\pp ^3}(m))=0$ pour tout $m$.
Donc la suite exacte
$$
0\rightarrow \Jj _{C/\pp ^3}(\sigma +n-\delta ) \rightarrow
\Jj _{C/\pp ^3}(\sigma +n) \rightarrow
\Jj _{P/X}(\sigma +n) \rightarrow 0,
$$
 donne l'implication:
$$
H^2(\Jj _{C/\pp ^3}(\sigma + n-\delta ))=0 \;  \Rightarrow \;
H^1(E'(n))=0
$$
(en effet,
d'apr\`es  la suite exacte $(\ast )$ on a
$h^1(E'(n)) \leq h^1(\Jj _{P/X}(\sigma + n))$).

De plus, la r\'esolution de  l'id\'eal homog\`ene de $C$ donne la
r\'esolution de
faisceaux
$$
0\rightarrow \Oo _{\pp ^3}(-s-1)^s \rightarrow  \Oo _{\pp ^3}(-s)^{s+1}
\rightarrow \Jj _{C/\pp ^3}\rightarrow 0,
$$
d'o\`u l'implication:
$$
H^3(\Oo _{\pp ^3}(\sigma + n -\delta -s-1))=0 \;  \Rightarrow \; H^1(E'(n))=0,
$$
ou encore:
$$
\sigma + n -\delta -s-1 \geq -3 \;  \Rightarrow \; H^1(E'(n))=0,
$$

Rappelons que $s= \delta - 2$ ou $\delta - 3$ et $\sigma = s/2$.
On en d\'eduit facilement que, pour $n> \frac{3}{2}\delta -4$,
on a $
\sigma + n -\delta -s-1 \geq -3$.
\eop

\begin{theoreme}
\label{A1}
Pour tout entier $c_2 > \frac{1}{12}(13 \delta ^3 - 24
\delta ^2 + 8 \delta )$,
le fibr\'e g\'en\'eral de la bonne  composante de
$M_X(c_2)$ que nous avons construite a la cohomologie naturelle.
\end{theoreme}
{\em D\'em:}
D'apr\`es la proposition \ref{R} on a $\beta \leq
\frac{3}{2}\delta - 3$.
Soit $P(t)$ le polyn\^ome tel que
\newline
$P(n) = 2\Cchi (\Oo _X(n))$.
Pour $t\geq k/2$, $P(t)$ est une fonction croissante  donc
 $$
\gamma := P(\beta ) \leq P( \frac{3}{2}\delta - 3).
$$
On conclut en calculant:
$$
P( \frac{3}{2}\delta - 3)=
\frac{1}{12}(13 \delta ^3 - 24 \delta ^2 + 8 \delta ).
$$
\eop

\end{document}